\documentclass[11pt]{article}
\usepackage{amsfonts,latexsym,rawfonts,amsmath,amssymb,amsthm,graphicx}
\numberwithin{equation}{section}
\newtheorem{theorem}{Theorem}[section]
\newtheorem{lem}[theorem]{Lemma}
\newtheorem{thm}[theorem]{Theorem}
\newtheorem{pro}[theorem]{Proposition}
\newtheorem{cor}[theorem]{Corollary}
\newtheorem{defi}[theorem]{Definition}
\newtheorem{rem}[theorem]{Remark}
\def\s{\,\,\,\,}

\def\mv{1.7ex}
\def\endproof{$\hfill\Box$\\}
\def\R{\mathbb{R}}

\def\div{{\rm div\,}}

\title{On Willmore Legendrian surfaces in $\mathbb{S}^5$ and the contact stationary Legendrian Willmore surfaces}
\author{Yong Luo}
\date{}
\begin{document}
\maketitle
\begin{abstract}In this paper we study  Willmore Legendrian surfaces (that is Legendrian surfaces which are critical points of the Willmore functional). We use an equality proved in \cite{Luo} to get a relation between Willmore Legendrian surfaces and contact stationary Legendrian surfaces in $\mathbb{S}^5$, and then we use this relation to prove a classification result for Willmore Legendrian spheres in $\mathbb{S}^5$. We also get an integral inequality for Willmore Legendrian surfaces and in particular we prove that if the square length of the second fundamental form of a Willmore Legendrian surface in $\mathbb{S}^5$ belongs to $[0,2]$, then it must either be $0$ and $L$ is totally geodesic or $2$ and $L$ is a flat minimal Legendrian tori, which generalizes a result of \cite{YKM}. We also study variation of the Willmore functional among Legendrian surfaces in 5-dimensional Sasakian manifolds. Let $\Sigma$ be a closed surface and $(M,\alpha,g_\alpha,J)$ a 5-dimensional Sasakian manifold with a contact form $\alpha$, an associated metric $g_\alpha$ and an almost complex structure $J$. Assume that $f:\Sigma\mapsto M$ is a Legendrian immersion. Then $f$ is called a contact stationary Legendrian Willmore surface (in short, a csL Willmore surface) if it is a critical point of the Willmore functional under contact deformations.  To investigate the existence of csL Willmore surfaces we introduce a higher order flow which preserves the Legendre condition and decreases the Willmore energy. As a first step we prove that this flow is well posed if $(M,\alpha,g_\alpha,J)$ is a Sasakian Einstein manifold, in particular $\mathbb{S}^5$.
\end{abstract}
\section{Introduction}
Let $\Sigma$ be a closed surface, $(M,g)$  a Riemannian manifold and $f$ an immersion from $\Sigma$ to $M$. We consider the Willmore functional defined as
\begin{eqnarray}
W(f)=\frac{1}{2}\int_\Sigma S-2|H|^2d\mu_f,
\end{eqnarray}
where $A$ is the second fundamental form of $f$ with respect to the induced metric, $S:=|A|^2$, $H$ is the mean curvature vector field of $f$ defined by
$$H=\frac{1}{2}tr A,$$
and $d\mu_f$ is the measure element on $f$.

We see that by the Gauss equation
$$4|H|^2-S=2K_f-2K_M,$$
where $K_f$ is the Gauss curvature of $f$ with respect to the induced metric and $K_M$ is the Gauss curvature of the metric $g$ when restricted to the plane $f_\ast(T\Sigma)$. Therefore by the Gauss-Bonnet theorem we see that $W(f)$ differs from the following \textbf{conformal Willmore functional}
\begin{eqnarray}
W_{conf}(f)=\int_\Sigma|H|^2+K_Md\mu_f
\end{eqnarray}
by a constant. It is well-known that $W_{conf}(.)$ is invariant under the conformal transformations of the ambient manifold $M$(cf.~\cite{Wh}, \cite{Ch}). Therefore $W(f)$ is also invariant under the conformal transformations of the ambient manifold $M$.

We see that $W(f)$ is nonnegative and it has the advantage of being zero exactly when $f$ is totally umbilical.

For a smooth variation $f: \Sigma\times I\rightarrow M$ and $\phi=\partial_tf$ we have (cf.~\cite{Tho}, \cite{Wei})
\begin{eqnarray}\label{First Variation}
\frac{d}{dt}W(f)=\frac{1}{2}\int_\Sigma(\langle \overrightarrow{W}(f),\phi\rangle) d\mu_f,
\end{eqnarray}
with $\overrightarrow{W}(f)=\sum_{\alpha=3}^n\overrightarrow{W}(f)^\alpha e_\alpha$, where $\{e_\alpha:3\leq\alpha\leq n\}$ is the orthonormal basis of the normal bundle of $f(\Sigma)$ in $M$ and
\begin{eqnarray}\label{Willmore operator}
\overrightarrow{W}(f)^\alpha&=&\Delta^{\nu}H^\alpha+\sum_{i,j,\beta}h_{ij}^\alpha h_{ij}^\beta H^\beta-2|H|^2H^\alpha, 3\leq\alpha\leq n,
\end{eqnarray}
where $\Delta^\nu$ is the Laplace-Beltrami operator along the normal vector bundle of $f$ and $h_{ij}^\alpha$ is the component of $A$ and $H^\alpha$ is a half of the trace of $(h_{ij}^\alpha)$.

A smooth immersion $f:\Sigma\rightarrow M$ is called a {\it Willmore immersion}, if it satisfies the {E-L} equation of $W$, i.e. $f$ is a Willmore immersion if and only if
\begin{eqnarray}\label{ELE}
\Delta^{\nu}H^\alpha+\sum_{i,j,\beta}h_{ij}^\alpha h_{ij}^\beta H^\beta-2|H|^2H^\alpha=0, 3\leq\alpha\leq n,
\end{eqnarray}
or equivalently,
\begin{eqnarray}
\Delta^\nu H+Q(A^\circ)H=0,
\end{eqnarray}
where $$Q(A^\circ)H:=\sum_{i,j,k,l}g^{ik}g^{jl}A^\circ_{ij}\langle A^\circ_{kl},H\rangle,$$ and $A^\circ:=A-Hg$ is the trace free part of $A$.

 In \cite{LW}, a geometrically constrained variation problem was studied for Lagrangian surfaces in $\mathbb{C}^2$. In this paper we introduce another kind of geometrically constrained variation problem of the Willmore functional, i.e. the contact deformations of the Willmore functional with $(M,\alpha,g_\alpha,J)$ being a 5-dimensional Sasakian manifold.
\begin{defi} A Legendrian immersion $f:\Sigma\mapsto (M,\alpha,g_\alpha,J)$ is called \textbf{a contact stationary Legendrian Willmore surface}(in short, a csL Willmore surface) if it is a critical point of the Willmore functional under contact deformations.
\end{defi}
In the following we will see that csL Willmore surfaces satisfy the following E-L equation:
\begin{eqnarray}\label{HSWS}
\div (J\overrightarrow{W}(f)-4JH)=0.
\end{eqnarray}
\begin{defi}
Let $(M^{2n+1}, \alpha, g_\alpha, J)$ be a Sasakian manifold. A Legendrian submanifold of  $M$ is called \textbf{a contact stationary Legendrian submanifold}(i.e. csL submanifold) if it is a critical point of the volume functional with respect to contact deformations.
\end{defi}
Contact stationary Legendrian submanifolds satisfy the following E-L equation (cf. \cite{CLU}\cite{Ir}):
\begin{eqnarray}
div_g{JH}=0£¬
\end{eqnarray}
where $g$ is the induced metric on $L$.
\begin{defi}
A Willmore and Legendrian surface in a 5-dimensional Sasakian manifold is called a Willmore Legendrian surface.
\end{defi}
In \cite{Luo} we proved that if $\Sigma$ is a Legendrian submanifold in a Sasakian manifold, then
\begin{eqnarray}\label{Main ob.}
\langle\Delta^\nu H, \textbf{R}\rangle=-2div_gJH,
\end{eqnarray}
where $\textbf{R}$ is the Reeb field. From (\ref{Main ob.}) we can prove in the following that
\begin{thm}
 Assume that $\Sigma$ is a Willmore Legendrian surface in a 5-dimensional Sasakian manifold $(M,\alpha,g_\alpha,J)$. Then it is a contact stationary Legendrian surface in $M$. In particular, every Willmore Legendrian surface in $\mathbb{S}^5$ with the standard Sasakian structure is a csL surface.
 \end{thm}

Based on this discovery, we can quite directly get in the following that
\begin{thm}[Theorem \ref{Legendrian sphere}]
Let $\Sigma$ be a Willmore Legendrian 2-sphere in the standard sphere $\mathbb{S}^5$. Then $\Sigma$ must be the equatorial 2-sphere.
\end{thm}
The classification of Willmore 2-spheres in $n$-sphere($n\geq3$) is a long-standing open problem, solved when $n = 3$ by Bryant(\cite{Br}), when $n=4$ by Ejiri(\cite{Ej}), Musso(\cite{Mu}) and Montiel(\cite{Mo}) independently, and recently the classification is solved when $n=5$ by Ma, Wang and Wang(\cite{MWW}). Castro and Urbano(\cite{CU}) proved that the Whitney sphere is the only Willmore Lagrangian sphere in $\mathbb{R}^4,$ which could be compared with our result.

In addition, following similar ideas from \cite{Luo} we can prove in the following that
\begin{thm}[Theorem \ref{Legendrian tori}]\label{thm2}
Let $\Sigma$ be a Willmore Legendrian surface in the standard sphere $\mathbb{S}^5$ and $0\leq S:=|A|^2\leq 2$, then it must either be totally geodesic or be a flat minimal Legendrian tori.
\end{thm}
\begin{rem}
Theorem  \ref{thm2} generalizes a related gap theorem for minimal Legendrian surfaces in $\mathbb{S}^5$ due to Yamaguchi, Kon and Miyahara(cf.~\cite{YKM}).
\end{rem}

The integral inequality and gap phenomenon for Willmore surfaces in $n$-sphere was first studied by Li(see \cite{Li1}-\cite{Li3}), where he generalized several integral inequalities and gap theorems for minimal surfaces in  $n$-sphere. In his papers Li used some Simons' type inequalities and the Willmore equation to get the desired integral inequalities and gap theorems. In our case by exploring the geometry of Legendrian surfaces we can get a new Simons' type inequality which helps us to find our integral inequality for Willmore Legendrian surfaces in $\mathbb{S}^5$.

In this paper we also introduce a flow which aims to prove the existence of csL Willmore surfaces.
 \begin{defi}
Let $f: \Sigma\times [0, T) \mapsto (M,\alpha,g_\alpha,J)$ be a family of immersions and $\textbf{R}$ the Reeb field of $M$. We call the following deformation the Legendrian Willmore flow(in short LeW flow)
 \begin{eqnarray}\left\{\begin{array}{rl}\label{flow11}
\frac{\partial f_t}{\partial t}&=-\div_t(J \overrightarrow{W}(f_t)-4JH)\textbf{R}-\frac{1}{2}J\nabla_t\div_t(J \overrightarrow{W}(f_t)-4JH),\\[\mv]
\\f(.,0)&=f_0,
\end{array}\right.
\end{eqnarray}
where $f(.,t)=f_t(.)$, $f_0$ is an initial Legendrian immersion, and $\div_t$, $\nabla_t$ are the divergence and gradient operators with respect to the induced metric on $f_t$ respectively.
\end{defi}
 This flow decreases the Willmore energy and preserves the Legendre condition. In the last section we will prove that it is equivalent to a sixth order scalar flow which is elliptic around the initial data and so by using a general existence theorem of Huisken and Polden(cf.~\cite{HP}) we can prove in the following the well posedness  of the LeW flow if the target manifold is a Sasakian Einstein manifold.
 \begin{thm}[Theorem \ref{main theorem}]
Let $(M,\alpha,g_\alpha)$ be a 5-dimensional Sasakian Einstein manifold. Then the LeW flow is well-posed, i.e. for any smooth Legendrian immersion $f_0$, there exists a $T>0$ and a unique family of Legendrian immersions $f_t$ for $t\in[0,T)$ such that $f_t$ satisfies (\ref{flow11}) with initial condition $f_0$.
\end{thm}
\begin{rem}
There is a very natural flow, that is the celebrated mean curvature flow, which preserves the Lagrange condition(see \cite{Sm1}). But it does not preserve the Legendre condition. In the past decade efforts were made by several authors to find flows preserving the Legendre condition. In \cite{Sm2} the author defined the Legendre mean curvature flow, which could be seem as a modification of the MCF and preserves the Legendre condition. In \cite{Le} a fourth order flow preserving the Legendre condition in spheres was defined. L$\hat{e}$'s work provides a general method to construct gradient flow preserving Legendre condtion and the definition of LeW flow is mainly motivated by her work. The LeW flow is another example whcih preserves the Legendre condtion for surfaces. The LeW flow is a Legendrian analogue of our recently defined HW flow in \cite{LW}.
\end{rem}

 \textbf{The rest of this paper is organized as follows:} In section 2 we give some basic material of contact geometry, which will be our geometry and analysis frame. In section 3 we give the definition of csL Willmore surfaces and prove a classification result of Willmore Legendrian spheres and  a gap theorem for Willmore Legendrian surfaces in $\mathbb{S}^5$. In section 4 we introduce the LeW flow which is closely related to the existence problem of csL Willmore surfaces and prove the well posedness of this flow when the target manifold is a Sasakian Einstein manifold.
\section{Basic material}
In this section we record some basic material of contact geometry. We invite the reader to consult \cite{Gei} and \cite{Bl} for more materials.
\subsection{Contact Manifolds}
\begin{defi}
A contact manifold $M$ is an odd dimensional manifold with a one form $\alpha$ such that $\alpha\wedge(d\alpha)^n\neq0$, where $dimM=2n+1$.
\end{defi}
Assume now that $(M,\alpha)$ is a given contact manifold of dimension $2n+1$. Then $\alpha$ defines a $2n-$dimensional vector bundle over $M$, where the fibre at each point $p\in M$ is given by
$$\xi_p=ker\alpha_p.$$
Sine $\alpha\wedge (d\alpha)^n$ defines a volume form on $M$, we see that
$$\omega:=d\alpha$$
is a closed nondegenerate 2-form on $\xi\oplus\xi$ and hence it defines a symplectic product on $\xi$, say $\omega$, such that $(\xi,\omega|_{\xi\oplus\xi})$ becomes a symplectic vector bundle. A consequence of this fact is that there exists an almost complex bundle structure $\tilde{J}:\xi\mapsto\xi$ compatible with $d\alpha$, i.e. a bundle endomorphism satisfying:
\\(1) $\tilde{J}^2=-id_\xi$,
\\(2) $d\alpha(\tilde{J}X,\tilde{J}Y)=d\alpha(X,Y)$ for all $X,Y\in\xi$,
\\(3) $d\alpha(X,\tilde{J}X)>0$ for $X\in\xi\setminus {0}$.

Since $M$ is an odd dimensional manifold, $\omega$ must be degenerate on $TM$, and so we obtain a line bundle $\eta$ over $M$ with fibres
$$\eta_p:=\{V\in TM|\omega(V,W)=0 \s\forall\s W\in\xi\}.$$
\begin{defi}
The Reeb vector field $\textbf{R}$ is the section of $\eta$ such that $\alpha(\textbf{R})=1$.
\end{defi}

Thus $\alpha$ defines a splitting of $TM$ into a line bundle $\eta$ with the canonical section $\textbf{R}$ and a symplectic vector bundle $(\xi,\omega|\xi\oplus\xi)$. We denote the projection along $\eta$ by $\pi$, i.e.
\begin{eqnarray*}
&&\pi:TM\mapsto\xi,
\\&&\pi(V):=V-\alpha(V)\textbf{R}.
\end{eqnarray*}
Using this projection we extend the almost complex structure $\tilde{J}$ to a section $J\in\Gamma(T^*M\otimes TM)$ by setting
$$J(V)=\tilde{J}(\pi(V)),$$
for $V\in TM$.

We have special interest in a kind of submanifolds in contact manifolds.
\begin{defi}
Let $(M,\alpha)$ be a contact manifold, a submanifold $\Sigma$ of $(M,\alpha)$ is called an isotropic submanifold if $T_x\Sigma\subseteq\xi$ for all $x\in \Sigma$.
\end{defi}
For algebraic reasons the dimension of an isotropic submanifold of a $2n+1$ dimensional contact manifold can not bigger than $n$.
\begin{defi}
An isotropic submanifold $\Sigma\subseteq(M,\alpha)$ of maximal possible dimension $n$ is called a Legendrian submanifold.
\end{defi}
\subsection{Sasakian manifolds}
Let $(M,\alpha)$ be a contact manifold with an almost complex structure $J$. A Riemannian metric $g_\alpha$ defined on $M$ is said to be associated, if it satisfies the following three conditions:
\\(1) $g_\alpha(\textbf{R},\textbf{R})=1$,
\\(2) $g_\alpha(V,\textbf{R})=0$, $\forall\s V\in\xi$,
\\(3) $\omega(V,JW)=g_\alpha(V,W)$, $\forall\s V,W\in\xi$.\

We should mention here that on any contact manifold there exists an associated metric on it, because we can construct one in the following way. We introduce a bilinear form $b$ by
$$b(V,W):=\omega(V,JW),$$
then the tensor
$$g:=b+\alpha\otimes\alpha$$
defines an associated metric on $M$.

Sasakian manifolds are the odd dimensional analogue of K\"ahler manifolds.
\begin{defi}
A contact manifold $M$ with an associated metric $g_\alpha$ is called Sasakian, if the cone $CM$ equipped with the following extended metric $\bar{g}$
\begin{eqnarray}\label{cone metic}
(CM,\bar{g})=(\mathbb{R}^+\times M,dr^2+r^2g_\alpha)
\end{eqnarray}
is K\"ahler w.r.t the following canonical almost complex structure $J$ on $TCM=\mathbb{R}\oplus\langle\textbf{R}\rangle\oplus\xi:$
$$J(r\partial r)=\textbf{R}, J(\textbf{R})=-r\partial r.$$
Furthermore if $g_\alpha$ is Einstein, $M$ is called a Sasakian Einstein manifold.
\end{defi}
\begin{lem}\label{cone cur.}
The Ricci curvature of the cone $(CM,\bar{g})$ satisfies (cf.~\cite{LeW}, A.3)
\begin{eqnarray*}
Ric(\partial r,\partial r)&=&-\frac{2n+1}{r}\frac{\partial^2r}{\partial r^2}=0,
\\Ric(\partial r,X)&=&0, \s if \s X\in T(\{r\}\times M),
\\Ric(X,Y)&=&Ric_M(\pi_\ast X,\pi_\ast Y)-2n\langle\pi_\ast X,\pi_\ast Y\rangle,\s if\s X,Y\in T(\{r\}\times M).
\end{eqnarray*}
\end{lem}
This shows if $(M,g_\alpha)$ is a Sasakian Einstein manifold then we can get a Calabi-Yau metric on $CM$ according to (\ref{cone metic}) by changing $g_\alpha$ to be some positive multiple of itself.

On Sasakian manifolds we have curvature conditions as follows (cf.~\cite{Bl}, Lemma 7.1 and Proposition 7.3, pp. 93-95).
\begin{lem}\label{curvature}
Let $(M,\alpha,g_\alpha,J)$ be a $2n+1$ dimensional Sasakian manifold. If $R$ and $Ric$ are the curvature tensor and the Ricci curvature of $(M,g_\alpha)$ respectively, then we have
\begin{eqnarray}
R(X,Y)\textbf{R}&=&\alpha(Y)X-\alpha(X)Y,
\\Ric(\textbf{R},X)&=&2m, \s if \s X=\textbf{R},
\\ Ric(\textbf{R},X)&=&0,\s if \s X\in ker\alpha,
\\and&&
\\R(X,Y)J Z &=&J(R(X,Y)Z)-g(Y,Z)JX \nonumber
\\&+&g(J X,Z)Y+g(X,Z)J Y-g(J Y,Z)X.
\end{eqnarray}
Moreover if $(M,g_\alpha)$ is Einstein, then the scalar curvature of $(M,g_\alpha)$ equals $2n(2n+1)$.
\end{lem}
We record more several lemmas which are well known in Sasakian geometry. These lemmas will be used in the subsequent sections.
\begin{lem}
Let $(M,\alpha,g_\alpha,J)$ be a Sasakian manifold. Then
\begin{eqnarray}\label{Reeb}
\bar{\nabla}_X\textbf{R}=-JX,
\end{eqnarray}
and
\begin{eqnarray}\label{derivatives}
(\bar{\nabla}_XJ)(Y)=g(X,Y)\textbf{R}-\alpha(Y)X,
\end{eqnarray}
for $X,Y\in TM$, where $\bar{\nabla}$ is the Levi-Civita connection on $(M,g_\alpha)$.
\end{lem}
\begin{lem}\label{mean curvatue form}
Let $\Sigma$ be a Legendrian submanifold in a Sasakian Einstein manifold $(M,\alpha,g_\alpha,J)$, then the mean curvature form $\omega(H,\cdot)|_\Sigma$ defines a closed one form on $\Sigma$.
\end{lem}
For a proof of this lemma we refer to \cite{Le}, Proposition A.2 and \cite{Sm2}, lemma 2.8. In fact they proved this result under the weaker assumption that $(M,\alpha,g_\alpha,J)$ is a weakly Sasakian Einstein manifold, where weakly Einstein means that $g_\alpha$ is Einstein only when restricted to the contact hyperplane.
\begin{lem}\label{orthogonal}
Let $\Sigma$ be a Legendrian submanifold in a Sasakian manifold $(M,\alpha,g_\alpha,J)$ and $A$ be the second fundamental form of $\Sigma$ in $M$. Then we have
\begin{eqnarray}
g_\alpha(A(X,Y),\textbf{R})=0.
\end{eqnarray}
\end{lem}
\proof For any $X,Y\in T\Sigma$,
\begin{eqnarray*}
\langle A(X,Y),\textbf{R}\rangle&=&\langle\bar{\nabla}_XY,\textbf{R}\rangle
\\&=&-\langle Y,\bar{\nabla}_X\textbf{R}\rangle
\\&=&\langle Y,JX\rangle
\\&=&0,
\end{eqnarray*}
where in the third equality we used (\ref{Reeb}).
\endproof

In particular this lemma implies that the mean curvature $H$ of $\Sigma$ is orthogonal to the Reeb field $\textbf{R}$. This fact is important in our following argument.
\begin{lem}\label{commute of J}
Let $(M,\alpha,g_\alpha,J)$ be a Sasakian manifold. For any $Y,Z\in Ker\alpha$, we have
\begin{eqnarray}
g_\alpha(\bar{\nabla}_X(JY),Z)=g_\alpha(J\bar{\nabla}_XY,Z).
\end{eqnarray}
\end{lem}
\proof Note that
$$(\bar{\nabla}_XJ)Y=\bar{\nabla}_X(JY)-J\bar{\nabla}_XY.$$
Therefore by using (\ref{derivatives}) we have
\begin{eqnarray*}
\langle \bar{\nabla}_X(JY),Z\rangle&=&\langle(\bar{\nabla}_XJ)Y,Z \rangle+\langle J\bar{\nabla}_XY,Z\rangle
\\&=&\langle J\bar{\nabla}_XY,Z\rangle,
\end{eqnarray*}
for any $Y,Z\in Ker\alpha$. \endproof

\textbf{The standard sphere $\mathbb{S}^{2n+1}$.}
 Let $\mathbb{C}^{n+1}=\mathbb{R}^{2n+2}$ be the Euclidean space with coordinates $(x_1,,...,x_{n+1},y_1,...,y_{n+1})$ and $\mathbb{S}^{2n+1}$ be the standard unit sphere in $\mathbb{R}^{2n+2}$. Define
$$\alpha_0=\sum_{j+1}^{n+1}(x_jdy_j-y_jdx_j),$$
then
$$\alpha:=\alpha_0|_{\mathbb{S}^{2n+1}}$$
defines a contact one form on $\mathbb{S}^{2n+1}$. Assume that $g_0$ is the standard metric on $\mathbb{R}^{2n+2}$ and $J_0$ is the standard complex structure of $\mathbb{C}^{n+1}$. We define
$g_\alpha=g_0|_{\mathbb{S}^{2n+1}},$
then $(\mathbb{S}^{2n+1},\alpha,g_\alpha)$ is a Sasakian Einstein manifold. The contact hyperplane is characterized by
$$ker\alpha_x=\{Y\in T_x\mathbb{S}^{2n+1}|\langle Y,J_0x\rangle=0\}.$$

\section{Willmore Legendrian surfaces and the Legendrian Willmore problem}
In this section we define the contact stationary Legendrian Willmore surfaces in 5-dimensional Sasakian manifolds.
\subsection{Definitions and the E-L equation.}
Let $(M,\alpha,g_\alpha,J)$ be a 5-dimensional Sasakian manifold with a contact structure $\alpha$, an associated metric $g_\alpha$, and an almost complex structure $J$. Assume that $\Sigma$ is a closed surface. An immersion $f$ from $\Sigma$ to $(M,\alpha,g_\alpha,J)$ is called a Legendrian immersion if $f(\Sigma)$ is a Legendrian surface of  $(M,\alpha,g_\alpha,J)$.
\begin{defi}
Let $(M,\alpha)$ be a contact manifold with contact 1-form $\alpha$ and $\Sigma$ be a Legendrian submanifold of $M$.
\\(1)A family $\{\Sigma_t:t\in[0,\epsilon),\epsilon>0\}$ of submanifolds in $M$ with $\Sigma_0=\Sigma$ is called a contact deformation of $\Sigma$ if $\Sigma_t$ are Legendrian submanifolds of $M$.
\\(2) A vector field $V$ on $\Sigma$ is called a contact variational vector field if there is a Legendrian deformation $\Sigma_t$ of $\Sigma$ such that $\frac{d\Sigma_t}{dt}|_{t=0}=V.$
\end{defi}
\begin{defi}
A Legendrian immersion $f$ from a closed surface $\Sigma$ to $(M,\alpha,g_\alpha,J)$ is called a contact stationary Legendrian Willmore immersion(in short, a csL Willmore immersion) if it is a critical point of the Willmore functional under contact deformations.
\end{defi}
By definition, an immersion $f$ is a contact stationary Legendrian Willmore immersion if and only if it satisfies
$$\int_\Sigma\langle \overrightarrow{W}(f), V\rangle d\mu_f=0,$$
for any contact variational vector field $V$ along $f$.

Recall that contact variational vector fields along $f$ are modeled by functions $s$ on $f$ in the following way(cf.~\cite{Le} and \cite{Sm2}):
\begin{eqnarray}
V=s\textbf{R}+\frac{1}{2}J\nabla s,
\end{eqnarray}
where $\textbf{R}$ is the Reeb vector field of $M$. We have

\begin{lem}\label{ort}
Let $\Sigma$ be a Legendrian submanifold of a Sasakian manifold $(M,\alpha,g_\alpha,J)$. Then
\begin{eqnarray}
\langle\overrightarrow{W}(f),\textbf{R}\rangle=-2\div_gJH,
\end{eqnarray}
where $\textbf{R}$ is the Reeb vector field on $M$ and $div_g$ is the divergence operator w.r.t the induced metric $g$ on $\Sigma$.
\end{lem}
\proof From the proof of Lemma 3.7 in \cite{Luo} we see that $\langle\Delta^\nu H, \textbf{R}\rangle=-2div JH$ and then by Lemma \ref{orthogonal} we complete the proof. \endproof

By the representation of $V$ and lemma \ref{ort} we can easily see that $f$ is a contact stationary Legendrian Willmore immersion if and only if
$$\int_\Sigma\div(J \overrightarrow{W}(f)-4JH)sd\mu_f=0,$$
for any function $s$ on $f$. Therefore we have
\begin{pro}
Let $\Sigma$ be a closed surface and $f: \Sigma\mapsto (M,\alpha,g_\alpha,J)$ be a Legendrian immersion into a Sasakian manifold $(M,\alpha,g_\alpha,J)$, then $f$ is a contact stationary Legendrian Willmore immersion if and only if
\begin{eqnarray}\label{E-L}
\div(J \overrightarrow{W}(f)-4JH)=0.
\end{eqnarray}
\end{pro}
Clearly every Willmore Legendrian immersion in a Sasakian manifold is a contact stationary Legendrian Willmore immersion, by definition. Among them we in particular have the minimal Legendrian hexagonal torus in $\mathbb{S}^5\subseteq \mathbb{C}^3$ as an example.

\textbf{The minimal Legendrian hexagonal torus.} The minimal Legendrian hexagonal torus is defined by
\begin{eqnarray*}
T=\{(\xi^1,\xi^2,\xi^3)\in \mathbb{C}^3: |\xi^i|^2=\frac{1}{3},i=1,2,3 \s and \s Im(\xi^1\xi^2\xi^3)=0\},
\end{eqnarray*}
 where $\xi^i$, i=1,2,3, are complex numbers. It is a  Willmore Legendrian surface in $\mathbb{S}^5$ and it is conjectured that this torus minimizes the Willmore energy among all Legendrian tori in $\mathbb{S}^5$.

\textbf{ Problem 1.} How to construct Willmore Legendrian surfaces in $\mathbb{S}^5$ which are not minimal Legendrian surfaces?

\textbf{Problem 2.} How to construct csL Willmore surfaces in $\mathbb{S}^5$ which are not Willmore Legendrian surfaces?
\subsection{Willmore Legendrian surfaces in $\mathbb{S}^5$}
In this section we give a classification result for Willmore Legendrian spheres and a gap theorem for Willmore Legendrian surfaces in $\mathbb{S}^5$.
\subsubsection{Willmore Legendrian spheres}
First we give a definition.
\begin{defi}
A contact stationary Legendrian submanifold $\Sigma$ in a Sasakian manifold \\$(M^{2n+1},\alpha,g_\alpha,J)$ is a stationary point of the volume functional under contact deformations.
\end{defi}
A contact stationary Legendrian submanifold $\Sigma$ satisfies the following E-L equation (cf. \cite{CLU}\cite{Ir}):
\begin{eqnarray}
\div_g(JH)=0,
\end{eqnarray}
where $g$ is the induced metric on $\Sigma$ and $div_g$ is the divergence operator.
\begin{defi}
A Willmore and Legendrian surface in a 5-dimensional Sasakian manifold is called a Willmore Legendrian surface.
\end{defi}
We have
\begin{thm}{\label{relation}}
Let $\Sigma$ be a Willmore Legendrian surface in a Sasakian 5-manifold $(M,\alpha,g_\alpha,J)$. Then it is a contact stationary Legendrian surface in $(M,\alpha,g_\alpha,J)$.
\end{thm}
\proof This is a direct consequence of Lemma \ref{ort}. \endproof

By this fact we have the following classification theorem:
\begin{thm}\label{Legendrian sphere}
 Let $\Sigma$ be a  Willmore Legendrian 2-sphere in the standard sphere $\mathbb{S}^5$. Then $\Sigma$ must be the equatorial 2-sphere.
\end{thm}
\proof By the last theorem we have that $\Sigma$ satisfies the following equation
\begin{eqnarray}
\div_g(JH)=0.
\end{eqnarray}
It is easy to see that this is equivalent to
\begin{eqnarray}\label{mean curvature1}
\delta(H\lrcorner\omega)=0,
\end{eqnarray}
where $\delta$ is the dual operator of $d$ on $\Sigma$ w.r.t the induced metric $g$.

On the other hand the mean curvature 1-form $H\lrcorner \omega$ of a Legendrian surface in a Sasakian Einstein manifold is closed, that is
\begin{eqnarray}\label{MC2}
d(H\lrcorner\omega)=0.
\end{eqnarray}
From (\ref{mean curvature1})-(\ref{MC2}),  we deduce that $H\lrcorner\omega$ is a harmonic 1-form on $\Sigma$. Since there is no non-trial harmonic 1-form on 2-sphere, we see that $H=0$, i.e. $\Sigma$ is a minimal Legendrian 2-sphere.  Therefore $\Sigma$ is the equatorial 2-sphere in $\mathbb{S}^5$, by Yau's result(cf.~\cite{Yau}).
\endproof
\subsubsection{A gap theorem}
Obviously minimal Legendrian surfaces in $\mathbb{S}^5$ are a special kind of Willmore Legendrian surfaces. For minimal Legendrian surfaces in $\mathbb{S}^5$, we have a gap theorem which states that any  minimal Legendrian surface with $0\leq S:=|A|^2\leq2$ must be $S=0$ or $S=2$(cf.~\cite{YKM}).  Here we generalize this result to Willmore Legendrian surfaces.
\begin{thm}\label{Legendrian tori}
Let $\Sigma$ be a Willmore Legendrian surface in $\mathbb{S}^5$ with $0\leq S:=|A|^2\leq2$, then either $S=0$ and $\Sigma$ is totally geodesic or $S=2$ and $\Sigma$ is a flat minimal Legendrian tori.
\end{thm}
\proof We postpone the proof of this theorem to the end of this paper. \endproof
\section{The LeW-flow}
In this section we introduce a flow method to investigate the existence of the contact stationary Legendrian Willmore immersions. As a first step we prove that this flow is well defined when the target manifold is a Sasakian Einstein manifold.
\begin{defi}
Let $f: \Sigma\times [0, T) \mapsto (M,\alpha,g_\alpha,J)$ be a family of immersions. We call $f$ a solution of the Legendrian Willmore flow (the LeW flow) if
\begin{eqnarray}
\left\{\begin{array}{rl}\label{flow}
\frac{\partial f_t}{\partial t}&=-\div_t(J \overrightarrow{W}(f_t)-4JH)\textbf{R}-\frac{1}{2}J\nabla_t\div_t(J \overrightarrow{W}(f_t)-4JH),\\[\mv]
\\f(.,0)&=f_0,
\end{array}\right.
\end{eqnarray}
where $f(.,t)=f_t(.)$, $f_0$ is an initial Legendrian immersion, and $\div_t$, $\nabla_t$ are the divergence and gradient operators with respect to the induced metric on $f_t$ respectively.
\end{defi}
The right hand side of the LeW flow defines a contact variation on $f_t$, as a consequence the LeW flow preserves the Legendrian condition.
\begin{pro}
The Willmore functional is decreasing along the LeW flow.
\end{pro}
\proof By (\ref{First Variation}) and lemma \ref{ort} we have
\begin{eqnarray*}
\frac{d W(f_t)}{dt}&=&\frac{1}{2}\int_\Sigma\langle\overrightarrow{W}(f_t),\frac{\partial f_t}{\partial t}\rangle d\mu_{f_t}
\\&=&\frac{1}{2}\int_\Sigma\langle\overrightarrow{W}(f_t), -\div(J \overrightarrow{W}(f_t)-4JH)\textbf{R}-\frac{1}{2}J\nabla\div(J \overrightarrow{W}(f_t)-4JH)\rangle d\mu_{f_t}
\\&=&-\frac{1}{4}\int_\Sigma|\div(J \overrightarrow{W}(f_t)-4JH)|^2d\mu_{f_t}.
\end{eqnarray*}
Hence the conclusion follows.
\endproof

In the following we introduce a $L^2$ metric defined by L\^{e}(cf.~\cite{Le}) on $\Lambda$, the set of all Legendrian surfaces in $(M,\alpha,g_\alpha,J)$, such that the Legendrian Willmore flow is a negative gradient flow with respect to this metric.

Let $\Sigma\in\Lambda$ and $V_1$, $V_2\in T_\Sigma\Lambda$, then we have
$$V_1=f_1\textbf{R}+\frac{1}{2}J\nabla f_1 \s and \s V_2=f_2\textbf{R}+\frac{1}{2}J\nabla f_2,$$
for some functions $f_1$, $f_2$ on $\Sigma$. We define
\begin{eqnarray}\label{L^2 metric}
\langle V_1, V_2 \rangle_\Lambda:=\frac{1}{4}\int_\Sigma  f_1f_2dvol_\Sigma.
\end{eqnarray}
Note that this $L^2$ metric is different from the usual one because it only takes into account the Reeb component.

We have
\begin{pro}
The gradient $\nabla W$ of the Willmore functional $W$ with respect to the $L^2$ metric defined by (\ref{L^2 metric}) is
$$\nabla W=h_f\textbf{R}+\frac{1}{2}\nabla h_f,$$
where
$$h_f=\div(J\overrightarrow{W}(f)-4JH).$$
\end{pro}
\proof Let $f: \Sigma\times I\mapsto (M,\alpha,g_\alpha,J)$ be a family of immersions and denote by $\Sigma_t=f_t(\Sigma)$ with
$$\frac{\partial \Sigma_t}{\partial t}|_{t=0}=h\textbf{R}+\frac{1}{2}J\nabla h:=V.$$
Then the gradient $\nabla W$ of the Willmore functional $W$ with respect to the $L^2$ metric $\langle. ,.\rangle_\Lambda$, by definition, is
\begin{eqnarray*}
\langle\nabla W, V\rangle_\Lambda&:=&\frac{\partial W(f_t)}{\partial t}|_{t=0}
\\&=&\frac{1}{2}\int_\Sigma\langle\overrightarrow{W}(f), h \textbf{R}+\frac{1}{2}J\nabla h\rangle
\\&=&\frac{1}{2}\int_\Sigma\langle\overrightarrow{W}(f), \frac{1}{2}J\nabla h\rangle-2h\div JH
\\&=&\frac{1}{4}\int_\Sigma h_fh:=\langle h_f\textbf{R}+\frac{1}{2}J\nabla h_f, V\rangle_\Lambda,
\end{eqnarray*}
where in the third equality we used lemma \ref{ort}. Thus $\nabla W=h_f\textbf{R}+\frac{1}{2}\nabla h_f.$
\endproof

Let $f_t$ be a solution of the Legendrian Willmore flow. Set $\Sigma_t=f_t(\Sigma)$, we have
\begin{eqnarray}\label{Gradient flow}
\frac{\partial \Sigma_t}{\partial t}=-\nabla W(\Sigma_t),
\end{eqnarray}
which implies that the Legendrian Willmore flow is a negative gradient flow with respect to the $L^2$ metric $\langle.,.\rangle_\Lambda$.

Let $J^1(\Sigma_0)$ be the standard contact manifold, i.e. the 1-jet bundle over $\Sigma_0$ with the canonical contact structure and $\varphi$ be a contactmorphism from an open neighborhood of the zero section $\Sigma_0(\subseteq U\subseteq J^1(\Sigma_0)$) to a small neighborhood $N_\epsilon(\Sigma_0)\subseteq M$. Then for $t$ small $\Sigma_t=\varphi(g_t, dg_t)$ for some function $g_t$ on $\Sigma_0$. Using this representation we have
\begin{lem}\label{equivalent}
Equation (\ref{Gradient flow}) locally (i.e. there exists $T>0$ such that for all $t\in[0, T)$) is equivalent to the equation
\begin{eqnarray}\label{scalar flow}
\frac{\partial g_t}{\partial t}=-\div(J\overrightarrow{W}(\varphi(g_t,dg_t))-4JH).
\end{eqnarray}
\end{lem}
\proof First we show that equation (\ref{Gradient flow}) implies equation (\ref{scalar flow}). For $t$ small $\Sigma_t$ will belongs to $N_\epsilon(\Sigma_0)$ and on $U$ we have a induced metric $\varphi^*g_\alpha$. We write down equation (\ref{Gradient flow}) as follows
\begin{eqnarray}\label{flow1}
\frac{\partial \varphi(g_t, dg_t)}{\partial t}=-\div(J\overrightarrow{W}(\varphi(g_t, dg_t))-4JH)\textbf{R}-\frac{1}{2}J\nabla\div(J\overrightarrow{W}(\varphi(g_t, dg_t))-4JH).
\end{eqnarray}
Let $\textbf{R}^1$ be the Reeb vector field of $J^1(\Sigma_0)$. The LHS of equation (\ref{flow1}) is the sum of the Reeb component $\varphi_*(\frac{\partial g_t}{\partial t}\textbf{R}^1)$ and the fiber component $\varphi_*(\frac{d}{dt}dg_t)$. The fiber component lies in the contact hyperplane in $J^1(\Sigma_0)$ and so it is orthogonal to the Reeb field $\textbf{R}^1$, with respect to the induced metric $\varphi^*g_\alpha$. So we have
$$\varphi_*(\frac{\partial g_t}{\partial t}\textbf{R}^1)=-\div(J\overrightarrow{W}(\varphi(g_t, dg_t))-4JH)\textbf{R}.$$
Noting that $\varphi_*\textbf{R}^1=\textbf{R}$, we have
$$\frac{\partial g_t}{\partial t}=-\div(J\overrightarrow{W}(\varphi(g_t, dg_t))-4JH),$$
which is equation (\ref{scalar flow}).

On the other hand, because $\frac{\partial dg_t}{\partial t}$ lies in the contact hyperplane, the Reeb component of $\frac{\partial(g_t, dg_t)}{\partial t}$ is $\frac{\partial g_t}{\partial t}\textbf{R}^1$. Therefore the Reeb component of $\frac{\partial \varphi(g_t, dg_t)}{\partial t}$ is $\varphi_*(\frac{\partial g_t}{\partial t}\textbf{R}^1)$, which is
$$-\div(J\overrightarrow{W}(\varphi(g_t, dg_t))-4JH)\textbf{R}.$$
So the Legendrian vector field  $\frac{\partial \varphi(g_t, dg_t)}{\partial t}$ must be
$$-\div(J\overrightarrow{W}(\varphi(g_t, dg_t))-4JH)\textbf{R}-\frac{1}{2}J\nabla\div(J\overrightarrow{W}(\varphi(g_t, dg_t))-4JH)$$
 and we complete the proof of the lemma.
\endproof

As a first step in the study of the LeW flow, we will prove the well posedness of the LeW flow as follows
\begin{thm}\label{main theorem}
Let $(M,\alpha,g_\alpha,J)$ be a 5-dimensional Sasakian Einstein manifold. Then the LeW flow is well-posed, i.e. for any smooth Legendrian immersion $f_0$, there exists a $T>0$ and a unique family of Legendrian immersions $f_t$ for $t\in[0,T)$ such that $f_t$ satisfies (\ref{flow}) with initial condition $f_0$.
\end{thm}
The standard odd dimensional unit sphere given in section 2 is a Sasakian Einstein manifold. Therefore we have
\begin{cor}
Let $(\mathbb{S}^5,\alpha,g_\alpha,J)$ be the standard contact sphere given in section 2. Then the LeW flow with $M=\mathbb{S}^5$ is well-posed.
\end{cor}

\textbf{Proof of Theorem \ref{main theorem}.} By lemma \ref{equivalent}, it suffices to prove that the flow (\ref{scalar flow}) exists and unique at a time interval $[0,T)$ for some $T>0$ with the initial condition $g_0=0$.

We can see that flow (\ref{scalar flow}) is a sixth-order quasilinear scalar flow. We shall use the following general existence theorem for higher order quasilinear scalar flow on compact manifolds, due to Huisken and Polden(cf.~\cite{HP}).
\begin{thm}[\cite{HP}, Theorem 7.15]
Suppose that for a smooth initial data $u_0$ the operator of $2p$ order
$$A(u)=A^{i_1j_1...i_pj_p}(x,u,\nabla u,...,\nabla^{2p-1}u)D_{i_1j_1...i_pj_p}$$
is smooth and strongly elliptic in a neighborhood of $u_0$. Then the evolution equation
$$D_tu=-A(u)u+b,$$
where $b=b(x,u,\nabla u,...,\nabla^{2p-1}u)$ is smooth, has a unique smooth solution on some interval $[0,T).$
\end{thm}
\begin{rem}
Here $A^{i_1j_1...i_pj_p}$ is strongly elliptic means that it is can be decomposed as
$$A^{i_1j_1...i_pj_p}=(-1)^pE^{i_1j_1}E^{i_2j_2}...E^{i_pj_p},$$
where the 2-form $E$ is strictly positive: $E\geq\lambda g$ for some $\lambda>0$.
\end{rem}
Since the evolution equation (\ref{scalar flow}) is a scalar quasilinear flow equation, in view of this theorem it suffices to show that it is parabolic around a neighborhood of $g_0=0$.
\begin{lem}
Denote the RHS of (\ref{scalar flow}) by $V_t=-\div(J\overrightarrow{W}(\varphi(g_t,dg_t)))$. Then
\begin{eqnarray}\label{ellipticity1}
V_t=-\Delta\div(JH)+\s lower\s order\s terms,
\end{eqnarray}
where $\Delta,\div,H$ are the Laplacian, divergence and mean curvature vector of $\varphi(g_t,dg_t)$ respectively.
\end{lem}
\proof It easy to see that $J\Delta^\nu H=\Delta JH$ and so by (\ref{Willmore operator}) we obtain
\begin{eqnarray}\label{1}
V_t=-\div\Delta JH+ lower\s order\s terms.
\end{eqnarray}
By the Ricci formula we have
\begin{eqnarray}\label{2}
\div\Delta JH=\Delta\div JH+lower\s order\s terms.
\end{eqnarray}
Inserting (\ref{2}) into (\ref{1}), we get (\ref{ellipticity1}).
\endproof

It is known that the mean curvature vector field $H(\Sigma)$ of a Lagrangian submanifold $\Sigma$ in a Calabi-Yau manifold $M^{2n}$ is symplectically dual to the angle form $d\theta$(cf.~\cite{HL}), i.e.
\begin{eqnarray}{\label{Lagangian angle}}
H=J(d\theta)^\sharp,
\end{eqnarray}
where $\sharp:\Omega^1(\Sigma)\mapsto Vect(\Sigma)$ is defined by
$$X(Y)=\langle X^\sharp,Y\rangle$$
and $\theta$ is the real part of the complex value of $Vol_C(TM^{2n})$. Here $Vol_C$ denotes the holomorphic complex volume form on $M^{2n}$.

If $\Sigma$ is a Lengendrian submanifold in a Sasakian manifold $(M,\alpha,g_\alpha,J)$, we have an analogue of (\ref{Lagangian angle}). More precisely we denote by $det(M)$ the determinant bundle of the contact plane bundle $ker\alpha$ over $M$ and by $\mathfrak{Leg}(M)$ the bundle of oriented Legendrian planes in $ker\alpha$. We also denote by det the following bundle map
\begin{eqnarray*}
det:\mathfrak{Leg}(M)&\mapsto&det(M)
\\\omega&\mapsto&\omega\wedge J\omega.
\end{eqnarray*}
We have
\begin{lem}[\cite{Le}]
Let us denote by $\tilde{\alpha}$ the canonical connection form on the determinant bundle $det(M)$ over a Sasakian manifold $(M,\alpha,g_\alpha,J)$. Then the mean curvature $H(\Sigma)$ of an oriented Legendrian submanifold $\Sigma\subseteq M$ is symplectically dual to
\begin{eqnarray}\label{Legendrian angle}
h_\Sigma=(det\circ\rho)^*\tilde{\alpha},
\end{eqnarray}
i.e. $h_\Sigma=J(H(\Sigma))^\sharp$. Here $\rho:\Sigma\mapsto\mathfrak{Leg}(M)$ denotes the Gauss map which sends each point $x\in \Sigma$ to the plane $T_x\Sigma$.
\end{lem}
 Set $\Sigma=\varphi(g_t,dg_t)$, we shall compute the symbol in an open simply connected domain on $\Sigma$. For the simplicity we shall denote this domain also by $\Sigma$.

Since $M$ is Sasakian Einstein, the form $h_\Sigma$ is closed, therefore the restriction of $detM$ to $\Sigma$ is a flat $\mathbb{S}^1-$bundle. Since $\Sigma$ is simply connected we can choose a trivialization
\begin{eqnarray}
\Pi:detM|_\Sigma\mapsto \mathbb{S}^1
\end{eqnarray}
which is compatible with this connection, i.e. $i^*\tilde{\alpha}=\Pi^*(d\theta)=d\Pi$, where $i$ denotes the embedding of $detM|_\Sigma$ to $detM$, and $d\theta$ is the canonical one form in the circle $\mathbb{S}^1$ with coordinate $\theta$. Thus we can rewrite (\ref{Legendrian angle}) as follows
\begin{eqnarray}
h_\Sigma=d(\Pi\circ det\circ\rho).
\end{eqnarray}
This equation implies that $H(\Sigma)=J\nabla(\Pi\circ det\circ\rho)$ and so we can rewrite equation (\ref{scalar flow}) as follows
\begin{eqnarray}\label{ellipticity2}
\frac{\partial g_t}{\partial t}=\Delta^2(\Pi\circ det\circ\rho)+\s lower \s order\s terms,
\end{eqnarray}
where in the above equality we used the result of lemma \ref{ellipticity1}.
\begin{lem}({\cite{Le}}, Lemma 5.7.)\label{ellipticity3}
The symbol of the linearization of $\Pi\circ det\circ\rho\circ\varphi(g_t)$ at $g_0=0$ is a positive multiple of the identity matrix.
\end{lem}
\proof Sketch: Because $g_\alpha$ defines a Sasakian Einstein  manifold, by lemma \ref{curvature} it has positive scalar curvature and so by lemma \ref{cone cur.} there exists a positive constant $\sigma$ such that $\sigma g_\alpha$ induces a Calabi-Yau metric on $CM$, and the new flow in the new metric $\sigma g_\alpha$ (\ref{scalar flow}) is a scaling of the old flow (\ref{scalar flow}).

Let us denote by $\Pi_0:detM\mapsto \mathbb{S}^1$ the canonical trivialization of $detM$ on the Calabi-Yau, and by $\Pi$ the trivialization of $det(C(M))|_{C(\Sigma)}$ which is induced from
$$\Pi:detM|_\Sigma\mapsto \mathbb{S}^1.$$
Since two trivializations are compatible with the canonical connection form $\alpha$ on $detCM$, so they are the same. Therefore the linearization $D(\Pi\circ det\circ\rho\circ\varphi)_0$ is equal to the restriction of the linearization $D(\Pi_0\circ i\circ\rho_1\circ C(\varphi))$ to homogenous functions, i.e. the set of functions $f(r,x)$ on $C(\varphi)$ with $f(r,x)=r^2f(x).$

First we compute the linearization of the angle function $\theta(L_f)$ on a Lagrangian submanifold $L$ in a Calabi-Yau manifold $N$, where $L_f$ is the deformation of $L$ by a function $f$ on $L$ via the following formula:
$$L_f=\phi(J\nabla f).$$
Here $J\nabla f$ is a Lagrangian submanifold in $NL$ and $\phi$ is a symplectomorphism $NL\mapsto N$ which equals to the identity on $L$. Then we get
\begin{eqnarray}
\theta(L_f)=\arg(e^{i\theta_L}\det(Id+\sqrt{-1}\nabla\nabla f)),
\end{eqnarray}
where $\nabla$ is the covariant derivative on $L$.

Noting that
$$\theta(L_{\epsilon f})=\arg(e^{i\theta_L}\det(Id+\sqrt{-1}\nabla\nabla \epsilon f))=\theta_L+\arg\det(Id+\sqrt{-1}\nabla\nabla \epsilon f),$$
and
$$\det(Id+\sqrt{-1}\nabla\nabla \epsilon f)=1+i\epsilon\Delta_Lf+o(\epsilon),$$
we have
$$\frac{d}{d\epsilon}\theta(L_{\epsilon f})|_{\epsilon=0}=\Delta_Lf.$$
Hence
\begin{eqnarray}
D(\Pi_0\circ i\circ\rho_1\circ C(\varphi)df)=\Delta_Lf.
\end{eqnarray}
In addition we get for $\tilde{f}=r^2f(x)$,
$$\Delta_{C(\Sigma)}\tilde{f}=\Delta_\Sigma f-2nf,$$
which proves our statement.
\endproof

By (\ref{ellipticity2}) and lemma \ref{ellipticity3} we see that the flow (\ref{scalar flow}) is a parabolic flow around the initial data and so we get theorem \ref{main theorem} directly from Huisken and Polden's theorem.

This finishes the proof of theorem \ref{main theorem}. \endproof

\section{Appendix: Proof of Theorem \ref{Legendrian tori}}
Let $\Sigma$ be a Legendrian surface in $\mathbb{S}^5$ with the induced metric $g$. Let $\{e_1,e_2\}$ be an orthogonal tangent frame on $\Sigma$ such that $\{e_1,e_2,Je_1,Je_2,\textbf{R}\}$ be an orthonormal frame on $\mathbb{S}^5$.

In the following we use indices $i,j,k,l,s,t,m$ and $\beta,\gamma$ such that
\begin{eqnarray*}
1\leq i,j,k,l,s,t,m&\leq&2,
\\1\leq\beta,\gamma&\leq&3,
\\ \gamma^\ast=\gamma+2,\s \beta^\ast&=&\beta+2.
\end{eqnarray*}

Let $A$ be the second fundamental form of $\Sigma$ in  $\mathbb{S}^5$ and define
\begin{eqnarray}
h_{ij}^k&=&g_\alpha(A(e_i,e_j),Je_k),
\\h^3_{ij}&=&g_\alpha(A(e_i,e_j),\textbf{R}).
\end{eqnarray}
Then
\begin{eqnarray}
h_{ij}^k&=&h_{ik}^j=h_{kj}^i,
\\h^3_{ij}&=&0.
\end{eqnarray}
The Gauss equations and Ricci equations are
\begin{eqnarray}
R_{ijkl}&=&(\delta_{ik}\delta_{jl}-\delta_{il}\delta_{jk})+\sum_s(h^s_{ik}h^s_{jl}-h^s_{il}h^s_{jk})\label{basic equation 1}
\\R_{ik}&=&\delta_{ik}+2\sum_sH^sh^s_{ik}-\sum_{s,j}h^s_{ij}h^s_{jk},
\\2K&=&2+4H^2-S,
\\R_{3412}&=&\sum_i(h_{i1}^1h_{i2}^2-h_{i2}^1h_{i1}^2)\nonumber
\\&=&\det h^1+\det h^2,
\end{eqnarray}
where $h^1,h^2$ are the second fundamental forms w.r.t. the directions $Je_1$, $Je_2$ respectively.

 In addition we have the following Codazzi equations and Ricci identities
\begin{eqnarray}
h^\beta_{ijk}&=&h^\beta_{ikj},
\\h^\beta_{ijkl}-h^\beta_{ijlk}&=&\sum_mh^\beta_{mj}R_{mikl}+\sum_mh^\beta_{mi}R_{mjkl}+\sum_\gamma h^\gamma_{ij}R_{\gamma^\ast\beta^\ast kl}.\label{basic equation 2}
\end{eqnarray}

Using these equations, we can get the following Simons' type inequality:
\begin{lem}\label{main lemma1}
Let $\Sigma$ be a Legendrian surface in $\mathbb{S}^5$. Then we have
\begin{eqnarray}\label{main lemma}
\frac{1}{2}\Delta\sum_{i,j,\beta}(h^\beta_{ij})^2&\geq&|\nabla^T h|^2-4|\nabla^\nu H|^2 +\sum_{i,j,k,\beta}(h^\beta_{ij}h^\beta_{kki})_j \nonumber
\\&+&S+2(1+H^2)\rho^2-\rho^4-\frac{1}{2}S^2,
\end{eqnarray}
where $|\nabla^T h|^2=\sum_{i,j,k,s}(h^s_{ijk})^2$ and $|\nabla^T H|^2=\sum_{i,s}(H^s_i)^2$.
\end{lem}
\proof Using equations from (\ref{basic equation 1}) to (\ref{basic equation 2}), we have
\begin{eqnarray}\label{simon type}
\frac{1}{2}\Delta\sum_{i,j,\beta}(h^\beta_{ij})^2
&=&\sum_{i,j,k,\beta}(h^\beta_{ijk})^2+\sum_{i,j,k,\beta}h^\beta_{ij}h^\beta_{kijk}\nonumber
\\&=&|\nabla h|^2-4|\nabla^\nu H|^2+\sum_{i,j,k,\beta}(h^\beta_{ij}h^\beta_{kki})_j+\sum_{i,j,l,k,\beta} h^\beta_{ij}(h^\beta_{lk}R_{lijk}+h^\beta_{il}R_{lj})\nonumber
\\&+&\sum_{i,j,k,\beta,\gamma} h^\beta_{ij}h^\gamma_{ki}R_{\gamma^\ast\beta^\ast jk}\nonumber
\\&=&|\nabla h|^2-4|\nabla^\nu H|^2+\sum_{i,j,k,s}(h^s_{ij}h^s_{kki})_j+2K\rho^2-2(\det h^1+\det h^2)^2\nonumber
\\&\geq&|\nabla h|^2-4|\nabla^\nu H|^2+\sum_{i,j,k,\beta}(h^\beta_{ij}h^\beta_{kki})_j+2(1+H^2)\rho^2-\rho^4-\frac{1}{2}S^2,
\end{eqnarray}
where $\rho^2:=S-2H^2$ and in the above calculations we used the following identities
\begin{eqnarray*}
\sum_{i,j,k,l,\beta} h^\beta_{ij}(h^\beta_{lk}R_{lijk}+h^\beta_{il}R_{lj})&=&2K\rho^2,
\\\sum_{i,j,k,\beta,\gamma} h^\beta_{ij}h^\gamma_{ki}R_{\gamma^\ast\beta^\ast jk}&=&-2(\det h^1+\det h^2)^2,
\end{eqnarray*}
where in the first equality we used $R_{lijk}=K(\delta_{lj}\delta_{ik}-\delta_{lk}\delta_{ij})$ and $R_{lj}=K\delta_{lj}$ in a proper coordinate, because $\Sigma$ is a surface.

Note that
\begin{eqnarray}\label{main idea1}
|\nabla h|^2&=&\sum_{i,j,k,\beta}(h^\beta_{ijk})^2\nonumber
\\&=&|\nabla^T h|^2+\sum_{i,j,k}(h^3_{ijk})^2\nonumber
\\&=&|\nabla^T h|^2+\sum_{i,j,k}(h^k_{ij})^2\nonumber
\\&=&|\nabla^T h|^2+S,
\end{eqnarray}
where in the third equality we used
\begin{eqnarray*}
h^3_{ijk}&=&\langle\bar{\nabla}_{e_k}A(e_i,e_j),\textbf{R}\rangle
\\&=&-\langle A(e_i,e_j),\bar{\nabla}_{e_k}\textbf{R}\rangle
\\&=&\langle A(e_i,e_j),Je_k\rangle
\\&=&h^k_{ij}.
\end{eqnarray*}
Combing (\ref{simon type}) and (\ref{main idea1}) we get (\ref{main lemma}). \endproof

Let $\Sigma$ be a surface in $\mathbb{S}^n$ with second fundamental form $A=(h_{ij}^\alpha), \alpha=3,...,n.$ We define the trace free tensor
$$\tilde{h}_{ij}^\alpha=h_{ij}^\alpha-H^\alpha g_{ij}.$$
Then the Willmore equation becomes
\begin{eqnarray}
\Delta^\nu H^\alpha+\sum_{\beta,i,j}\tilde{h}_{ij}^\alpha\tilde{h}_{ij}^\beta H^\beta=0, 3\leq\alpha\leq n.
\end{eqnarray}
We have
\begin{lem}\label{main lemma2}
Let $\Sigma$ be a Willmore surface in $\mathbb{S}^n$, then
\begin{eqnarray}\label{equ1}
\int_\Sigma|\nabla^\nu H|^2d\mu=\int_\Sigma\tilde{\sigma}_{\alpha\beta}H^\alpha H^\beta d\mu,
\end{eqnarray}
where $\tilde{\sigma}_{\alpha\beta}=\sum_{ij}\tilde{h}_{ij}^\alpha\tilde{h}_{ij}^\beta$.
\end{lem}
\proof See \cite{Li2}. \endproof

Because $(\tilde{\sigma}_{\alpha\beta})$ is a symmetric matrix we can assume that it is diagonal at a point $p\in \Sigma$, by choosing a proper local orthonormal frame around $p$. Then we see that $\tilde{\sigma}_{\alpha\beta}=\tilde{\sigma}_\alpha\delta_{\alpha\beta},\rho^2=\sum_\alpha\tilde{\sigma}_\alpha$ and
\begin{eqnarray}\label{ine1}
H^2\rho^2=\sum_\alpha (H^\alpha)^2\sum_\beta\tilde{\sigma}_\beta\geq\sum_\alpha (H^\alpha)^2\tilde{\sigma}_\alpha=\sum_{\alpha,\beta}H^\alpha H^\beta\tilde{\sigma}_{\alpha\beta}.
\end{eqnarray}
Integrating over (\ref{main lemma}) and using equality (\ref{equ1}), inequality (\ref{ine1}) we get
\begin{eqnarray*}
0&=&\int_\Sigma|\nabla^T h|^2-4|\nabla^\nu H|^2+S+2(1+H^2)\rho^2-\rho^4-\frac{1}{2}S^2d\mu
\\ &\geq& \int_\Sigma-4|\nabla^\nu H|^2+S+2(1+H^2)\rho^2-\rho^4-\frac{1}{2}S^2d\mu
\\&=& \int_\Sigma-4|\nabla^\nu H|^2+4H^2\rho^2+\rho^2(2-S)+S-\frac{1}{2}S^2d\mu
\\&\geq&\int_\Sigma(\rho^2+\frac{S}{2})(2-S)d\mu.
\end{eqnarray*}
Therefore if $0\leq S\leq 2$ we must have $S=0$, i.e. $\Sigma$ is totally geodesic or $S=2$.

At last we analyze the case $S=2$. In this case we must have
$$(\det h^1+\det h^2)^2=\frac{S^2}{4},$$
that's $\det h^1+\det h^2=\frac{-S}{2}$ or $\frac{S}{2}$. If $\det h^1+\det h^2=\frac{S}{2}$, noting that $\rho^2=S-2H^2=\frac{S}{2}-(\det h^1+\det h^2)=0$, which implies $|\nabla^\nu H|^2=0$, by (\ref{equ1}) and (\ref{ine1}). But similar to (\ref{main idea1}) we have $|\nabla^\nu H|^2=|\nabla^TH|^2+H^2$, we have $H=0$ and so $S=0$, a contradiction. Thus we must have $\det h^1+\det h^2=\frac{-S}{2}$. Noting that $\rho^2=S-2H^2=\frac{S}{2}-(\det h^1+\det h^2)=S$ and $\rho^2=S-2H^2$, we get $H=0$, i.e. $\Sigma$ is a flat minimal Legendrian tori. This completes the proof of Theorem \ref{Legendrian tori}. \endproof
\quad\\

\textbf{Acknowledgement.} The author started this project when he was a Ph.D. student of professor Guofang Wang at Albert-Ludwigs Universit\"at Freiburg. He is very appreciated with Guofang Wang for stimulating discussions and constant support. Many thanks to the referee for his/her comments and suggestions which made this paper more readable. The author is partially supported by the NSF of China(No.11501421).

{}
\vspace{1cm}\sc
Yong Luo

School of Mathematics and statistics,

Wuhan University, Wuhan 430072, China

{\tt yongluo@whu.edu.cn}

\vspace{1cm}\sc


\begin{thebibliography}{2}\small
\bibitem[Bl]{Bl} D. E. Blair, Reimannian geometry of contact and symplectic manifolds, {\em Progress in Math., vol 203,} Birkh\"auser, Basel, 2002.
\bibitem[Br]{Br} R. Bryant, A duality theorem for Willmore surfaces, {\em J. Diff. Geom.} {\bf20}(1984), 23--53.
\bibitem[CLU]{CLU} I. Castro, H. Z. Li and F. Urbano, Hamiltonian-minimal lagrangian submanfolds in complex space form, {\em Pacific J. Math.} {\bf227}(2006), 43-65.
\bibitem[CU]{CU} I. Castro and F. Urbano, Willmore Surfaces of $\mathbb{R}^4$ and the Whitney Sphere, {\em Ann. Global Anal. Geom.} {\bf19}(2001), 153--175.
\bibitem[Ch]{Ch} B. Y. Chen, Some conformal invariants of submanifolds and their applications, {\em Boll. Un. Mat. Ital.} {\bf(4) 10}(1974), 380--385.
\bibitem[Ej]{Ej} N. Ejiri, Willmore surfaces with a duality in $\mathbb{S}^N(1)$, {\em Proc. Math. Lond. Soc.} {\bf(3)57}(1988), 383--416.
\bibitem[Gei]{Gei}H. Geiges, An introduction to contact topology, {\em Cambridge studies in advanced mathematics,} vol. {\bf109}, Cambridge University Press, Cambridge, 2008.

\bibitem[HL]{HL}R. Harvey, B. Lawson, Calibrated geometries, {\em Acta Math.} {\bf148}(1982), 48--157.
\bibitem[HP]{HP} G. Huisken, A. Polden, Geometric evolution equations for hypersurfaces, in: {\em lecture notes in math.,} {Vol. \bf 1713} (1999), 45--84.
\bibitem[Ir]{Ir} H. Iriyeh, Hamitonian minimal Lagrangian cones in $\mathbb{C}^m$, {\em Tokyo J. Math.} {\bf28}(2005), 91--107.
\bibitem[Le]{Le} H. V. L\^e, A minimizing deformation of Legendrian submanifolds in the standard sphere, {\em Diff. Geom. and its Appl.} {\bf 21}(2004), 297--316.
\bibitem[LeW]{LeW}H. V. L\^e, and G. F. Wang, Anti-complexified Ricci flow on compact symplectic manifolds, {\em J. Reine Angew. Math.} {\bf530}(2001), 17--31.
\bibitem[Li1]{Li1} H. Z. Li, Willmore hypersurfaces in a sphere, {\em Asian J. Math.} {\bf5}(2001), 365--378.
\bibitem[Li2]{Li2} H. Z. Li, Willmore Surfaces in $\mathbb{S}^n$, {\em Ann. Global Anal. Geom.} {\bf8}(2002), 203--213.
\bibitem[Li3]{Li3} H. Z. Li, Willmore submanifolds in a sphere, {\em Mathematical Research Letters} {\bf9}(2002), 771--790.
\bibitem[Luo]{Luo} Y. Luo, Contact stationary Legendrian surfaces in $\mathbb{S}^5$, arXiv:1211.4227v6, submitted.
\bibitem[LW]{LW}Y. Luo and G. F. Wang, On geometrically constrained variational problems of the Willmore functional I: the Lagrangian Willmore problem, {\em Comm. Anal. Geom.} {\bf23}(2015), 191--223.
\bibitem[MWW]{MWW} X. Ma, C. P. Wang and P. Wang, Classification of Willmore 2-spheres in the 5-dimensional sphere, arXiv:1409.2427v2, preprint(2014).
\bibitem[Mo]{Mo} S. Montiel, Willmore two-spheres in the four-sphere, {\em Trans. Amer. Math. Soc.} {\bf352}(2000), 4469--4486.
\bibitem[Mu]{Mu} E. Musso, Willmore surfaces in the four-sphere, {\em Ann. Global Anal. Geom.} {\bf8}, No.1(1990), 21--41.
\bibitem[Sm1]{Sm1}  K. Smoczyk, Angle theorems for the Lagrangian mean curvature flow, {\em Math. Z.} {\bf240}(2002), 849--883.
\bibitem[Sm2]{Sm2} K. Smoczyk, Closed Legendre geodesics in Sasaki manifolds, {\em New York J. Math.} {\bf9}(2003), 23--47.
\bibitem[Tho]{Tho} G. Thomsen, \"Uber konforme Geometrie I ; Grundlagen der konformen Flachentheorie, {\em Ab. Math. Sem. Univ. Hamburg} {\bf3}(1924), no. 1, 31--56.
\bibitem[Wei]{Wei} J. Weiner, On a problem of Chen, Willmore, et al., {\em Indiana Univ. Math. J.} {\bf 27}(1978), 19--35.
\bibitem[Wh]{Wh} J. H. White, A global invirant of conformal mappings in space, {\em Proc. Amer. Math. Soc.} {\bf38}(1973), 162--164.
 \bibitem[YKM]{YKM} S. Yamaguchi, M. Kon and Y. Miyahara, A theorem on C-totally real minimal surface, {\em Proc. Amer. Math. Soc.}{\bf54}(1976), 276-280.
 \bibitem[Yau]{Yau} S. T. Yau, Submanifolds with constant mean curvature I, {\em Amer. J. of Math.} {\bf96}(1974), 346--366.
\end{thebibliography}
\end{document}